# Infinitely ramified Galois representations

By Ravi Ramakrishna

In this paper we show how to construct, for most $p \geq 5$, two types of surjective representations $\rho : G_\mathbb{Q} = \text{Gal}(\bar{\mathbb{Q}}/\mathbb{Q}) \to \text{GL}_2(\mathbb{Z}_p)$ that are ramified at an infinite number of primes. The image of inertia at almost all of these primes will be torsion-free. The first construction is unconditional. The catch is that we cannot say whether $\rho \mid_{G_p = \text{Gal}(\bar{\mathbb{Q}}_p/\mathbb{Q}_p)}$ is crystalline or even potentially semistable. The second construction assumes the Generalized Riemann Hypothesis (GRH). With this assumption we can further arrange that $\rho \mid_{G_p}$ is crystalline at $p$. We remark that infinitely ramified *reducible* representations have been previously constructed by more elementary means.

We outline the method. Let $E_{/\mathbb{Q}}$ be a (modular!) semistable elliptic curve with good reduction at 3. Let $p > 3$ be a prime of good ordinary reduction such that for all $l$ prime, $v_l(j(E))$ is *not* divisible by $p$ where $j(E)$ is the $j$-invariant of $E$. Assume also that $a_p \neq \pm 1$. The collection of such $p$ form a set of density 1 (see [M2]). Let $S_0$ be the set containing all primes of bad reduction of $E$, $p$, and the infinite prime. For a set $T$ of primes denote by $G_T$ the Galois group over $\mathbb{Q}$ of the maximal extension of $\mathbb{Q}$ unramified outside places of $T$.

Suppose the residual representation $\bar{\rho} : G_{S_0} \to \text{GL}_2(\mathbb{F}_p)$ arising from the Galois action on the $p$-torsion of $E$ is surjective. (Since $E$ does not have complex multiplications Serre has shown in [Se2] this is the case for almost all $p$.) The set of $p$ satisfying all the above conditions is density 1. Let $\text{Ad}\,\bar{\rho}$ be the set of $2 \times 2$ matrices in $\mathbb{F}_p$ where Galois acts through $\bar{\rho}$ and by conjugation. Recall the exact sequence

$$0 \to \text{III}^2_{S_0}(\text{Ad}\,\bar{\rho}) \to H^2(G_{S_0}, \text{Ad}\,\bar{\rho}) \to \oplus_{v \in S_0} H^2(G_v, \text{Ad}\,\bar{\rho}).$$

In [Fl], Flach gives a condition, holding for all but finitely many $p$, that guarantees that $\text{III}^2_{S_0}(\text{Ad}\,\bar{\rho})$ is trivial. Mazur has shown in [M2] (for our chosen $p$) that $H^2(G_v, \text{Ad}\,\bar{\rho}) = 0$ for all $v \in S_0$. Thus for $p$ in a set of density 1 we have that $H^2(G_{S_0}, \text{Ad}\,\bar{\rho})$ is trivial. This is significant as obstructions to lifting problems lie in $H^2(G_{S_0}, \text{Ad}\,\bar{\rho})$. Henceforth assume $p \geq 5$ is a prime satisfying all of these conditions.

Let $\rho_0$ be the Galois representation associated to the $p$-adic Tate module of the elliptic curve $T_p(E)$. We then inductively construct a sequence $\{\rho_k\}$ of surjective representations of $G_\mathbb{Q}$ onto $\mathrm{GL}_2(\mathbb{Z}_p)$, with $\rho_{k-1} \equiv \rho_k \bmod p^k$ and $\rho_k$ ramified at a new prime $l_k$ (or possibly two new primes $l_{k1}$ and $l_{k2}$ in the GRH case). This ramification at the new prime(s) in $\rho_k$ will first appear mod $p^{k+1}$. For every $k$ we ensure that $\det\rho_k = \chi$, the cyclotomic character. For our purposes this restriction means it suffices to study the cohomology of $\mathrm{Ad}^0\bar\rho$, the $2 \times 2$ matrices in $\mathbb{F}_p$ with trace zero, as opposed to that of $\mathrm{Ad}\,\bar\rho$. That the image of inertia at $l_k$ (or $\{l_{k1}, l_{k2}\}$) is infinite follows from the fact that $\mathrm{GL}_2(\mathbb{Z}_p)$ has no torsion elements congruent to $I$ mod $p$ for $p > 2$. (For $l \neq p$, the pro-$p$ part of the inertia group at $l$ is pro-cyclic.)

We will also arrange in our GRH result for each $\rho_k$ to be ordinary at $p$. In this construction $\rho_k \mid_{G_p} = \begin{pmatrix} \psi\chi & * \\ 0 & \psi^{-1} \end{pmatrix}$ for all $k \geq 0$ where $\chi$ is the cyclotomic character and $\psi$ is unramified with $\psi^2 \neq 1$. Up to isomorphism over $\mathbb{Q}_p$ there is only one such (local at $p$) nontrivial representation and it is crystalline. Alternatively, by the theorems of [W] and [TW], $\rho_k$ is modular of weight 2 and level prime to $p$ and therefore crystalline at $p$.

In both the unconditional and GRH cases the limit of the $\rho_k$ will be our $\rho$. (In the GRH case the limit will be the ordinary representation above and therefore crystalline at $p$.) The main theorems are stated below with Theorem 2 in a slightly simplified form.

THEOREM 1. *Fix $E_{/\mathbb{Q}}$ a semistable elliptic curve with good reduction at 3. For primes $p \geq 5$ in a set of density one, there exist surjective representations $G_\mathbb{Q} \to \mathrm{GL}_2(\mathbb{Z}_p)$ ramified at infinitely many primes. The reduction mod $p$ of these representations is the Galois action on the $p$-torsion of $E$.*

THEOREM 2. *Assume the GRH. Consider $E$ as above. For primes $p \geq 5$ in a set of density one there exist surjective representations $G_\mathbb{Q} \to \mathrm{GL}_2(\mathbb{Z}_p)$ ramified at infinitely many primes that are crystalline at $p$. The reduction mod $p$ of these representations is again the Galois action on the $p$-torsion of $E$.*

*Acknowledgement.* I would like to thank the referee for several helpful suggestions and Jim Cogdell for advice and encouragement.

*Deformation theory.* We give a short introduction to deformation theory. See [M1], [M3], [BM], [B1] and [B2] for details and more results.

Let $\bar\pi : H \to \mathrm{GL}_d(\mathbb{F}_q)$ be an absolutely irreducible continuous representation of a profinite group $H$ where $\mathbb{F}_q$ is the finite field of $q$ elements. Suppose $H^1(H, \mathrm{Ad}\,\bar\pi)$ is finite-dimensional. Let $\mathcal{C}$ be the category of Artinian local rings with residue field $\mathbb{F}_q$ where the morphisms are homomorphisms that induce the identity map on the residue field. Let $R$ be in $\mathcal{C}$. We call two lifts $\gamma_1$ and $\gamma_2$ of $\bar{\pi}$ to $\mathrm{GL}_n(R)$ strictly equivalent if $\gamma_1 = A\gamma_2 A^{-1}$ for some $A$ congruent to the identity matrix modulo the maximal ideal $m_R$ of $R$. We call a strict equivalence class of lifts of $\bar{\pi}$ to $R$ a deformation of $\bar{\pi}$ to $R$.

Mazur studied the deformations of $\bar{\pi}$ and proved the following fundamental theorem in [M1].

THEOREM A. *There is a complete local Noetherian ring $R^{\mathrm{un}}$ with residue field $\mathbb{F}_q$ and a continuous homomorphism $\tilde{\pi} : H \to \mathrm{GL}_d(R^{\mathrm{un}})$ such that:*

1. *Reduction of $\tilde{\pi}$ modulo the maximal ideal of $R^{\mathrm{un}}$ gives $\bar{\pi}$.*

2. *For any ring $R$ in $\mathcal{C}$ and any deformation $\gamma$ of $\bar{\pi}$ to $\mathrm{GL}_n(R)$ there is a unique homomorphism $\phi : R^{\mathrm{un}} \to R$ in $\mathcal{C}$ such that $\phi \circ \tilde{\pi} = \gamma$ as deformations.*

Moreover, if $\bar{\pi}$ is not absolutely irreducible the statements hold except that the $\phi$ in part 2 may not be unique. We call $R^{\mathrm{un}}$ the universal deformation ring associated to $H$ and $\bar{\pi}$ in the absolutely irreducible case. We call $R^{\mathrm{un}}$ the versal ring associated to $H$ and $\bar{\pi}$ otherwise.

Let $W(\mathbb{F}_q)$ be the ring of Witt vectors of $\mathbb{F}_q$. In either case we have the following fact.

*Fact.* $R^{\mathrm{un}}$ is a quotient of $W(\mathbb{F}_q)[[T_1, T_2, ...T_r]]$ where
$$r = \dim_{\mathbb{F}_q} H^1(H, \mathrm{Ad}\,\bar{\pi}).$$

The elements of $H^1(H, \mathrm{Ad}\,\bar{\pi})$ correspond to the deformations of $\bar{\pi}$ to $\mathbb{F}_q[\varepsilon] = \mathbb{F}_q[X]/(X^2)$, the dual numbers of $\mathbb{F}_q$. Given $f \in H^1(H, \mathrm{Ad}\,\bar{\pi})$ the corresponding lift to the dual numbers is given by $\pi_f(\sigma) = (I + \varepsilon f(\sigma))\bar{\pi}(\sigma)$.

We now specialize the situation. Assume $\mathbb{F}_q = \mathbb{F}_p$ for some prime $p$ and that our representations are two-dimensional. Let $\pi_n$ be a deformation of $\bar{\pi}$ to $\mathrm{GL}_2(\mathbb{Z}/p^n)$. We may ask whether $\pi_n$ deforms to $\mathrm{GL}_2(\mathbb{Z}/p^{n+1})$. The obstruction to deforming $\pi_n$ to $\mathrm{GL}_2(\mathbb{Z}/p^{n+1})$ lies in $H^2(H, \mathrm{Ad}\,\bar{\pi})$. If this obstruction is trivial $\pi_n$ deforms to some $\pi_{n+1}$ and $\mathrm{pr} \circ \pi_{n+1} = \pi_n$ where $\mathrm{pr} : \mathbb{Z}/p^{n+1} \to \mathbb{Z}/p^n$ is the canonical projection. In the unobstructed case one sees that $H^1(H, \mathrm{Ad}\,\bar{\pi})$ acts on the set of deformations of $\pi_n$ to $\mathrm{GL}_2(\mathbb{Z}/p^{n+1})$. For $f \in H^1(H, \mathrm{Ad}\,\bar{\pi})$ the action is given by $(f.\pi_{n+1})(\sigma) = (I + p^n f(\sigma))(\pi_{n+1}(\sigma))$. If $\bar{\pi}$ is absolutely irreducible $H^1(H, \mathrm{Ad}\,\bar{\pi})$ acts on the the deformations of $\pi_n$ to $\mathrm{GL}_2(\mathbb{Z}/p^{n+1})$ as a principal homogeneous space.

Mazur also showed that modifications could be made so that related functors with the *ordinary* restriction were also representable. Here $H$ is a Galois group and we insist that when restricted to a suitable inertia group $I$ we only consider lifts $\pi$ of $\bar{\pi}$ whose restriction to $I$ is of the form $\begin{pmatrix} \psi & * \\ 0 & 1 \end{pmatrix}$. See [M1] for details.

*Local at $l$ deformation theory.* Let $l \notin S_0$ be a prime (at which we eventually wish to allow ramification) and let $G_l = \text{Gal}(\bar{\mathbb{Q}}_l/\mathbb{Q}_l)$. Suppose $\bar{\rho} : G_{\mathbb{Q}} \to \text{GL}_2(\mathbb{F}_p)$ is unramified at $l$ and $\bar{\rho}\,|_{G_l}$ is given by $\bar{\rho}(\sigma_l) = \begin{pmatrix} 2 & 0 \\ 0 & 1 \end{pmatrix}$ where $\sigma_l$ corresponds to Frobenius at $l$. Since we assume that the determinant is the cyclotomic character we have $l \equiv 2 \bmod p$. (Our choice of 2 is arbitrary. Any value $\neq \pm 1$ will serve our purposes. This is one reason why we insist $p \neq 3$.)

LEMMA 1. $H^2(G_l, \text{Ad}^0\bar{\rho})$ *is one-dimensional and* $H^1(G_l, \text{Ad}^0\bar{\rho})$ *is two-dimensional.*

*Proof.* Note that with $G_l$ action, $\text{Ad}^0\bar{\rho} \simeq \mathbb{F}_p \oplus \mu_p \oplus \mu_p(-1)$. As $l \equiv 2 \bmod p$, we see $\mu_p$ are not contained in $\mathbb{Q}_l$ and $H^0(G_l, \text{Ad}^0\bar{\rho})$ is one-dimensional.

By local duality we see $H^0(G_l, \text{Ad}^0\bar{\rho}^*)$ and $H^2(G_l, \text{Ad}^0\bar{\rho})$ are dual where $X^*$ is by definition $\text{Hom}(X, \mu_p)$ with Galois action. Since $p \geq 5$, $\mu_p(-1)$ and $\mu_p$ are not isomorphic as $G_l$ modules so we see $H^0(G_l, \text{Ad}^0\bar{\rho}^*)$ is one-dimensional. An application of the local Euler characteristic gives the result for $H^1(G_l, \text{Ad}^0\bar{\rho})$.

We want to consider deformations of $\bar{\rho}$ to $\mathbb{Z}/p^n$. As $l \neq p$, such deformations factor through the Galois group of the maximal tamely ramified extension of $\mathbb{Q}_l$ over $\mathbb{Q}_l$. This group is well understood. (See [Se1].) Thus we may assume that $G_l$ is topologically generated by $\sigma_l$ and $\tau_l$ subject to the relation $\sigma_l \tau_l \sigma_l^{-1} = \tau_l^l$ where $\tau_l$ topologically generates inertia and, as above, $\sigma_l$ corresponds to Frobenius.

*Definition.* We say $\rho : G_l \to \text{GL}_2(\mathbb{Z}/p^n)$ is *special* if $\rho$ is given by $\sigma_l \mapsto \begin{pmatrix} l & 0 \\ 0 & 1 \end{pmatrix}$ and $\tau_l \mapsto \begin{pmatrix} 1 & u \\ 0 & 1 \end{pmatrix}$ for $u \in \mathbb{Z}/p^n$.

*Remark.* In practice $u$ will be a nonzero multiple of $p$. In previous papers we used the expression *desired form* in similar situations. Here we use the term *special* for the sake of consistency with the terminology of modular forms.

The images of $\sigma_l$ and $\tau_l$ satisfy the relation above. Note the problem of deforming $\rho$ to $\mathbb{Z}/p^{n+1}$ is obviously unobstructed. One need only lift $u$ from mod $p^n$ to mod $p^{n+1}$ to get a special deformation of $\rho$ to mod $p^{n+1}$.

We give a basis for $H^1(G_l, \text{Ad}^0\bar{\rho})$. Recall that

$$\bar{\rho}(\sigma_l) = \begin{pmatrix} 2 & 0 \\ 0 & 1 \end{pmatrix} \text{ and } \bar{\rho}(\tau_l) = \begin{pmatrix} 1 & 0 \\ 0 & 1 \end{pmatrix}.$$

A nontrivial unramified lift to $\mathbb{F}_p[\varepsilon]$ is given by

$$\rho(\sigma_l) = \begin{pmatrix} 2 & 0 \\ 0 & 1 \end{pmatrix} + \varepsilon \begin{pmatrix} 2 & 0 \\ 0 & -1 \end{pmatrix}, \ \rho(\tau_l) = I + \varepsilon \begin{pmatrix} 0 & 0 \\ 0 & 0 \end{pmatrix}.$$

This corresponds to the unramified 1-cohomology class given by

$$r_l(\sigma_l) = \begin{pmatrix} 1 & 0 \\ 0 & -1 \end{pmatrix}, \ r_l(\tau_l) = \begin{pmatrix} 0 & 0 \\ 0 & 0 \end{pmatrix}.$$

A nontrivial ramified lift to the dual numbers is given by

$$\rho(\sigma_l) = \begin{pmatrix} 2 & 0 \\ 0 & 1 \end{pmatrix} + \varepsilon \begin{pmatrix} 0 & 0 \\ 0 & 0 \end{pmatrix}, \ \rho(\tau_l) = I + \varepsilon \begin{pmatrix} 0 & 1 \\ 0 & 0 \end{pmatrix}.$$

Since $l \equiv 2 \mod p$ we see that $\sigma_l \tau_l \sigma_l^{-1} = \tau_l^l$ holds. The corresponding 1-cohomology class is given by

$$s_l(\sigma_l) = \begin{pmatrix} 0 & 0 \\ 0 & 0 \end{pmatrix}, \ s_l(\tau_l) = \begin{pmatrix} 0 & 1 \\ 0 & 0 \end{pmatrix}.$$

Note that both $r_l$ and $s_l$, or more precisely their corresponding deformations to $\mathbb{F}_p[\varepsilon]$, cut out $\mathbb{Z}/p$ extensions of $\mathbb{Q}_l(\bar{\rho})$, the extension of $\mathbb{Q}_l$ fixed by the kernel of $\bar{\rho}\mid_{G_l}$. Any nontrivial linear combination of $r_l$ and $s_l$ cuts out the unique $\mathbb{Z}/p \times \mathbb{Z}/p$ extension of $\mathbb{Q}_l(\bar{\rho})$.

Also note that for any special lift of $\bar{\rho}\mid_{G_l}$ to mod $p^n$, $n \geq 2$, acting on it by the 1-cohomology class $s_l$ preserves specialness. One sees this by noting

$$(I + p^{n-1}s(\sigma_l))\begin{pmatrix} l & 0 \\ 0 & 1 \end{pmatrix} = \begin{pmatrix} l & 0 \\ 0 & 1 \end{pmatrix}$$

and

$$(I + p^{n-1}s(\tau_l))\begin{pmatrix} 1 & u \\ 0 & 1 \end{pmatrix} = \begin{pmatrix} 1 & u + p^{n-1} \\ 0 & 1 \end{pmatrix}.$$

Thus acting on a special local at $l$ deformation by a multiple of $s_l$ leaves the local at $l$ lifting problem unobstructed. For this reason we call $s_l$ a *null* 1-cohomology class. If $n - 1 > k$ and $u \neq 0$ then acting on a deformation to mod $p^n$ by $s_l$ gives a new deformation that is still ramified at $l$.

PROPOSITION 1. *Let $\bar{\rho}\mid_{G_l}$ be unramified and given by $\bar{\rho}(\sigma_l) = \begin{pmatrix} 2 & 0 \\ 0 & 1 \end{pmatrix}$. Fix $f \in H^1(G_l, \mathrm{Ad}^0 \bar{\rho})$ independent of the null 1-cohomology class $s_l$. Let $\rho_n$ be a special (at $l$) deformation of $\bar{\rho}\mid_{G_l}$ to mod $p^n$ and $\rho_{n+1}\mid_{G_l}$ be any (local at $l$) deformation of $\rho_n$ to mod $p^{n+1}$. Then there is an $\alpha \in \mathbb{F}_p$ such that $(\alpha f).\rho_{n+1}$ is special at $l$ (and itself unobstructed). Thus $\bar{\rho}\mid_{G_l}$ can be deformed to $\mathbb{Z}_p$ one step at a time with adjustments made at each step only by a multiple of $f$.*

*Proof.* Note that $\rho_{n+1}$ differs from a special deformation of $\rho_n$ to $\mathbb{Z}/p^{n+1}$ by the action of *some* element of the two-dimensional space $H^1(G_l, \mathrm{Ad}^0\bar{\rho})$. Since the one-dimensional subspace of null 1-cohomology classes preserves specialness we need only alter by a multiple of a nonnull 1-cohomology class, namely $f$. It is possible that in our characteristic zero representation we may have $\tau_l \mapsto I$, that is it might be unramified.

*Global considerations.* Recall that $\bar{\rho} : G_\mathbb{Q} \to \mathrm{GL}_2(\mathbb{F}_p)$ satisfies the numerous hypotheses of the introduction. The construction is inductive. Let $S_n = S_{n-1} \cup \{l_n\}$ where $l_n \equiv 2 \bmod p$ is as in the previous section and satisfies other conditions described later. The fact below follows immediately from our hypotheses on $\bar{\rho}$, Proposition 1.6 of [W], triviality of $\mathrm{III}^2_{S_0}(\mathrm{Ad}^0\bar{\rho})$ and Lemma 1.

*Fact* 1. The image of $H^1(G_{S_k}, \mathrm{Ad}^0\bar{\rho})$ in $H^1(G_{S_{k+1}}, \mathrm{Ad}^0\bar{\rho})$ under the (*injective*) *inflation map is codimension* 1 *and* $H^1(G_{S_n}, \mathrm{Ad}^0\bar{\rho})$ *is* $(n+2)$-*dimensional.*

Suppose for $0 \le n \le k$ we have constructed a surjective representation $\rho_n : G_{S_n} \to \mathrm{GL}_2(\mathbb{Z}_p)$ and that $\rho_{n-1} \equiv \rho_n \bmod p^n$. Suppose further that for $1 \le i \le n$, $\rho_n \mid_{G_{l_i}}$ is given by

$$\sigma_{l_i} \mapsto \begin{pmatrix} l_i & 0 \\ 0 & 1 \end{pmatrix} \text{ and } \tau_{l_i} \mapsto \begin{pmatrix} 1 & p^i u_{i,n} \\ 0 & 1 \end{pmatrix}$$

with $u_{i,n} \in \mathbb{Z}_p^*$. Also assume that for each $n$, $0 \le n \le k$ there exists $f_n \in H^1(G_{S_{n-1}}, \mathrm{Ad}^0\bar{\rho})$ that does not inflate from $H^1(G_{S_{n-2}}, \mathrm{Ad}^0\bar{\rho})$ such that $f_n \mid_{G_{l_n}}$ is nonnull. We also require that for $0 \le i < j \le n$ that $f_i \mid_{G_j}$ is trivial.

Our aim is to construct $l_{k+1}$, $\rho_{k+1}$ with $\rho_k \equiv \rho_{k+1} \bmod p^{k+1}$ and $f_{k+1} \in H^1(G_{S_k}, \mathrm{Ad}^0\bar{\rho})$ with $f_{k+1} \mid_{G_{l_{k+1}}}$ nonnull. We will show for $1 \le n \le k$ that $f_n \mid_{G_{l_{k+1}}}$ is trivial. This allows us to continue the induction. Then the limit $\rho$ of the $\{\rho_n\}$ exists.

Clearly $\rho_k$ factors through $G_{S_k}$. Let $\mathbb{Q}(\bar{\rho})$ denote the extension of $\mathbb{Q}$ cut out by the $p$-torsion of $E$ and $\mathbb{Q}(\rho_{k,k+2})$ the field cut out by $\rho_k \bmod p^{k+2}$. Let $\mathbb{Q}(\bar{\rho}_{k,\varepsilon})$ be the composite field cut out by all lifts to the dual numbers $\mathbb{F}_p[\varepsilon]$ that factor through $G_{S_k}$. Note that $\mathbb{Q}(\bar{\rho}_{k,\varepsilon})$ is closely linked to $H^1(G_{S_k}, \mathrm{Ad}^0\bar{\rho})$. Let $\mathbb{K}_k$ be the composite $\mathbb{Q}(\rho_{k,k+2})\mathbb{Q}(\bar{\rho}_{k,\varepsilon})$ and

$$C_k = \mathrm{Gal}(\mathbb{Q}(\rho_{k,k+2})/\mathbb{Q}) \simeq \mathrm{GL}_2(\mathbb{Z}/p^{k+2}), \ N_k = \mathrm{Gal}(\mathbb{Q}(\bar{\rho}_{k,\varepsilon})/\mathbb{Q}(\bar{\rho})).$$

LEMMA 2. *If* $p \ge 5$ *then* $\mathbb{Q}(\rho_{k,k+2}) \cap \mathbb{Q}(\bar{\rho}_{k,\varepsilon}) = \mathbb{Q}(\bar{\rho})$.

*Proof.* Let $\mathbb{L}$ denote this intersection. Suppose the intersection $\mathbb{L}$ strictly contains $\mathbb{Q}(\bar{\rho})$. Since $\mathbb{Q}(\rho_{k,k+2})$ and $\mathbb{Q}(\bar{\rho}_{k,\varepsilon})$ are both Galois over $\mathbb{Q}$, so is $\mathbb{L}$. Since $\bar{\rho}$ is onto $\mathrm{GL}_2(\mathbb{F}_p)$ we easily see $\mathrm{Ad}^0\bar{\rho}$ is irreducible as a $\mathrm{Gal}(\mathbb{Q}(\bar{\rho})/\mathbb{Q})$ module. As deformations of $\bar{\rho}$ to $\mathbb{F}_p[\varepsilon]$ factor through $\mathrm{Gal}(\mathbb{Q}(\bar{\rho}_{k,\varepsilon})/\mathbb{Q})$ we see the composition series for $N_k = \operatorname{Gal}(\mathbb{Q}(\bar{\rho}_{k,\varepsilon})/\mathbb{Q}(\bar{\rho}))$ as a $\operatorname{Gal}(\mathbb{Q}(\bar{\rho})/\mathbb{Q})$ module consists entirely of $\operatorname{Ad}^0\bar{\rho}$'s. Thus there is a field $\mathbb{J}$ between $\mathbb{L}$ and $\mathbb{Q}(\bar{\rho})$ with $\operatorname{Gal}(\mathbb{J}/\mathbb{Q}(\bar{\rho})) \simeq \operatorname{Ad}^0\bar{\rho}$ as $\operatorname{Gal}(\mathbb{Q}(\bar{\rho})/\mathbb{Q})$-modules. As $\mathbb{J} \subseteq \mathbb{Q}(\bar{\rho}_{k,\varepsilon})$, the sequence

$$1 \to \operatorname{Gal}(\mathbb{J}/\mathbb{Q}(\bar{\rho})) \to \operatorname{Gal}(\mathbb{J}/\mathbb{Q}) \to \operatorname{Gal}(\mathbb{Q}(\bar{\rho})/\mathbb{Q}) \to 1$$

corresponds to an element of $H^1(G_{S_k}, \operatorname{Ad}^0\bar{\rho})$ and thus splits. But $\mathbb{J} \subseteq \mathbb{Q}(\rho_{k,k+2})$ and thus $\operatorname{Gal}(\mathbb{J}/\mathbb{Q})$ is a quotient of $C_k = \operatorname{GL}_2(\mathbb{Z}/p^{k+2})$. Thus we see that $\operatorname{Gal}(\mathbb{J}/\mathbb{Q}) \simeq \operatorname{GL}_2(\mathbb{Z}/p^2)$. For $p \geq 5$ it is a simple exercise to see this is a nonsplit extension. This contradiction proves the lemma.

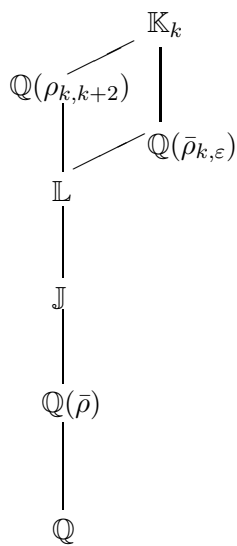

The diagram is therefore as below.

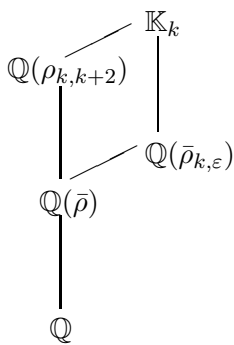

LEMMA 3. $\operatorname{Gal}(\mathbb{K}_k/\mathbb{Q}) \simeq$ *the semidirect product of* $C_k$ *by* $N_k$.

*Proof.* We abuse notation and denote $\text{Gal}(\mathbb{K}_k/\mathbb{Q}(\rho_{k,k+2}))$ by $N_k$. This is isomorphic to $\text{Gal}(\mathbb{Q}(\bar\rho_{k,\varepsilon})/\mathbb{Q}(\bar\rho))$ by Lemma 2. Denote $\text{Gal}(\mathbb{K}_k/\mathbb{Q}(\bar\rho_{k,\varepsilon}))$ by $M_k$ and $\text{Gal}(\mathbb{K}_k/\mathbb{Q})$ by $H_k$. By Lemma 2 we have the exact sequence

$$1 \to N_k \times M_k \to H_k \to \text{GL}_2(\mathbb{F}_p) = \text{Gal}(\mathbb{Q}(\bar\rho)/\mathbb{Q}) \to 1.$$

We also have the split exact sequence

$$1 \to (N_k \times M_k)/M_k \to H_k/M_k \to \text{GL}_2(\mathbb{F}_p) = \text{Gal}(\mathbb{Q}(\bar\rho)/\mathbb{Q}) \to 1.$$

The easiest way to see that the last sequence splits is to consider the universal deformation ring $R^{\text{un},k}$ associated to this problem and its maximal ideal $m_{R^{\text{un},k}}$. The deformation to the characteristic $p$ ring $R^{\text{un},k}/(p, m^2_{R^{\text{un},k}})$ factors through $H_k/M_k = \text{Gal}(\mathbb{Q}(\bar\rho_{k,\varepsilon})/\mathbb{Q})$ and $\text{GL}_2(\mathbb{F}_p) = \text{Gal}(\mathbb{Q}(\bar\rho)/\mathbb{Q})$ embeds in $\text{GL}_2(R^{\text{un},k}/(p, m^2_{R^{\text{un},k}}))$. Let $D$ be an image of $\text{GL}_2(\mathbb{F}_p)$ in $H_k/M_k$ associated to a splitting and $\tilde D$ the corresponding subgroup of $H_k$. We claim that in the exact sequence

$$1 \to N_k \to H_k \to H_k/N_k = C_k \to 1$$

the subgroup $\tilde D \subseteq H_k$ maps isomorphically to $C_k$. This will give the desired splitting. Counting orders it suffices to show $\tilde D \cap N_k = \{1\}$. But since $D \cap (N_k \times M_k)/M_k = \{M_k/M_k\}$ we are done.

For $a$ nonzero in $\mathbb{Z}/p$ let $a^*$ denote the Teichmüller lift of $a$ to $\mathbb{Z}_p$, i.e. the unique $p-1$st root of unity in $\mathbb{Z}_p$ congruent to $a$ mod $p$. For $p \geq 5$ we see $2^* \neq (1/2)^*$. Let $A \in \text{GL}_2(\mathbb{Z}/p^{k+2})$ be the matrix $\begin{pmatrix} 2^* & 0 \\ 0 & 1 \end{pmatrix}$. Let $B = \begin{pmatrix} 1 & 0 \\ 0 & -1 \end{pmatrix} \in \text{Ad}^0\bar\rho$ be an element of $N_k$ whose projection to $N_{k-1}$ is trivial. We are using that $\text{Gal}((\mathbb{Q}(\bar\rho)/\mathbb{Q}) \simeq \text{GL}_2(\mathbb{F}_p)$ acts on $N_k$. Such a $B$ exists because $N_{k-1}$ is a quotient of $N_k$ by a $\text{Gal}(\mathbb{Q}(\bar\rho)/\mathbb{Q})$ stable subgroup of $N_k$ isomorphic to $\text{Ad}^0\bar\rho$. This follows from Fact 1. Recall that $\text{Gal}(\mathbb{Q}(\bar\rho)/\mathbb{Q})$ acts on $N_k$ via $\bar\rho$ and conjugation. (When $k=0$ let $B = \begin{pmatrix} 1 & 0 \\ 0 & -1 \end{pmatrix}$ be *any* such element of $N_0$.) An application of Chebotarev's theorem gives the lemma below.

LEMMA 4. *If $p \geq 5$ then there are infinitely many primes whose Frobenius in $\text{Gal}(\mathbb{K}_k/\mathbb{Q})$ is in the conjugacy class of $(A, B)$ in the semidirect product $H_k$ of $C_k$ by $N_k$. Such elements have order $pc$ where $c$, which is prime to $p$, is the order of $2^*$ in $\mathbb{Z}_p^*$.*

*Remark.* It is important for our purposes that $(A, B)$ is *not* conjugate in $H_k$ to some $\left(\tilde A, \begin{pmatrix} 0 & 0 \\ 0 & 0 \end{pmatrix}\right)$. This is because we will be interested in how primes in our Chebotarev class split from $\mathbb{Q}$ to $\mathbb{Q}(\bar{\rho})$ to $\mathbb{Q}(\bar{\rho}_{k,\varepsilon})$. We do not want these primes to split completely from $\mathbb{Q}(\bar{\rho})$ to $\mathbb{Q}(\bar{\rho}_{k,\varepsilon})$. Our choice of $(A, B)$ guarantees this nonconjugacy and the corresponding nonsplitting.

The prime we wish to "add to the level" will be as in Lemma 4. We may choose *any* prime in our Chebotarev class as $l_{k+1}$. Since $\det \rho_k = \chi$, the cyclotomic character, we see such primes are congruent to $2^*$ mod $p^{k+2}$. We have the exact sequence

$$0 \to \text{III}^2_{S_{k+1}}(\text{Ad}^0 \bar{\rho}) \to H^2(G_{S_{k+1}}, \text{Ad}^0 \bar{\rho}) \to \oplus_{v \in S_{k+1}} H^2(G_v, \text{Ad}^0 \bar{\rho}).$$

As $\text{III}^2_{S_0}(\text{Ad}^0 \bar{\rho})$ is trivial (by assumption), $S_0 \subseteq S_{k+1}$ and $\text{III}^2_{S_0}(\text{Ad}^0 \bar{\rho}) \to \text{III}^2_{S_{k+1}}(\text{Ad}^0 \bar{\rho})$ is surjective by global Tate duality, the $\text{III}^2_{S_{k+1}}(\text{Ad}^0 \bar{\rho})$ term is trivial. Thus we need only analyze local lifting problems to analyze global lifting problems. As $p$ is odd there will be no obstructions at the Archimedean prime. Since for $v \in S_0$ we are assuming that $H^2(G_v, \text{Ad}^0 \bar{\rho})$ is trivial we only study primes in $S_{k+1} - S_0$.

We want $\rho_{k+1} |_{G_{l_{k+1}}}$ to be given by

$$\sigma_{l_{k+1}} \mapsto \begin{pmatrix} l_{k+1} & 0 \\ 0 & 1 \end{pmatrix}, \quad \tau_{l_{k+1}} \mapsto \begin{pmatrix} 1 & p^{k+1} u_{k+1,k+1} \\ 0 & 1 \end{pmatrix}$$

for some $u_{k+1,k+1} \in \mathbb{Z}_p^*$. Here $\tau_{l_{k+1}}$ topologically generates inertia at $l_{k+1}$ and $\sigma_{l_{k+1}}$ corresponds to Frobenius. Since $\sigma_{l_{k+1}}$ and $\tau_{l_{k+1}}$ satisfy the well-known relation $\sigma_{l_{k+1}} \tau_{l_{k+1}} \sigma_{l_{k+1}}^{-1} = \tau_{l_{k+1}}^{l_{k+1}}$ so must their images.

We start with $\rho_k$ mod $p^{k+2}$. We must adjust matters at this first stage so ramification actually occurs at $l_{k+1}$. This is done by altering $\rho_k$ mod $p^{k+2}$ by *any* nonzero element of $H^1(G_{S_{k+1}}, \text{Ad}^0 \bar{\rho})$ that does not inflate from $H^1(G_{S_k}, \text{Ad}^0 \bar{\rho})$. Such an element is clearly ramified at $l_{k+1}$. This guarantees that $\tau_{l_{k+1}}$ has nontrivial image as we lift to mod $p^m$ for all $m \geq k+3$. We need only do this once. Secondly, we need a global 1-cohomology class $f_{k+1} \in H^1(G_{S_{k+1}}, \text{Ad}^0 \bar{\rho})$ that is locally at $l_{k+1}$ independent of the null 1-cohomology class for $l_{k+1}$. That is, we must have a global 1-cohomology class that locally at $l_{k+1}$ satisfies the hypotheses of Proposition 1. Then we can alter our representation by a suitable multiple of this global 1-cohomology class to make the local representation at $G_{l_{k+1}}$ special. This procedure must be done as we lift from mod $p^m$ to mod $p^{m+1}$ for every $m \geq k+2$. We will also need $f_i |_{G_{l_{k+1}}}$ trivial for $1 \leq i \leq k$.

In [R2] we used the same 1-cohomology class for these two tasks. It was ramified at the prime in question, and brute computer computation showed it to be independent of the null 1-cohomology class locally. (Note that in [R2] we worked with $p = 3$ and there were various technical differences from our current situation.) Here, for each prime, we use different 1-cohomology classes for these two tasks. If for each prime $l_i$ we could use the same global 1-cohomology class for both purposes we could unconditionally construct odd representations ramified at an infinite number of primes that were crystalline at $p$. For the second task, the global 1-cohomology class will be unramified at $l_k$. This works because by Proposition 1 the null 1-cohomology class $s_{l_k}$ is ramified at the prime in question and thus independent of an unramified 1-cohomology class.

We perform the induction with $\rho_0$ as our starting point. We assume the existence of $\{\rho_0, \rho_1, \ldots, \rho_k\}$ with $\rho_n : G_{S_n} \to \mathrm{GL}_2(\mathbb{Z}_p)$ surjective and $\rho_{n-1} \equiv \rho_n$ mod $p^n$ for $0 \leq n \leq k$. Furthermore we insist that for $1 \leq i \leq n \leq k$, $\rho_n \mid_{G_{l_i}}$ is given by $\sigma_{l_i} \mapsto \begin{pmatrix} l_i & 0 \\ 0 & 1 \end{pmatrix}$ and $\tau_{l_i} \mapsto \begin{pmatrix} 1 & p^i u_{i,n} \\ 0 & 1 \end{pmatrix}$ with $u_{i,n} \in \mathbb{Z}_p^*$. Also for each $n$ with $1 \leq n \leq k$ we assume there exists $f_n \in H^1(G_{S_{n-1}}, \mathrm{Ad}^0 \bar{\rho})$ such that $f_n \mid_{G_{l_n}}$ is independent of the null 1-cohomology class at $l_n$ and that $f_j \mid_{G_{l_n}}$ is trivial for $j < n$.

To complete the induction we must construct $\rho_{k+1}$ with
$$\rho_k \equiv \rho_{k+1} \bmod p^{k+1},$$
$\rho_{k+1} \mid_{G_{l_i}}$ given as above for $1 \leq i \leq k$. We must also find $f_{k+1} \in H^1(G_{S_k}, \mathrm{Ad}^0 \bar{\rho})$ with $f_{k+1} \mid_{G_{l_{k+1}}}$ independent of the null 1-cohomology class at $l_{k+1}$ and $f_i \mid_{G_{l_{k+1}}}$ trivial for $i < k+1$.

LEMMA 5. *Let $l_{k+1}$ be any prime with Frobenius as in Lemma* 4. *There exists $f_{k+1} \in H^1(G_{S_k}, \mathrm{Ad}^0 \bar{\rho})$ such that $f_{k+1} \mid_{G_{l_{k+1}}}$ is independent of the null* 1-*cohomology class in $H^1(G_{l_{k+1}}, \mathrm{Ad}^0 \bar{\rho})$.*

*Proof.* Every $g \in H^1(G_{S_k}, \mathrm{Ad}^0 \bar{\rho})$ is clearly unramified at $l_{k+1}$. Suppose for all such $g$ we have that $g \mid_{G_{l_{k+1}}}$ is trivial. Then the primes above $l_{k+1}$ in $\mathbb{Q}(\bar{\rho})$ split completely in the composite field cut out by the deformations to the dual numbers factoring through $G_{S_k}$, that is they split completely from $\mathbb{Q}(\bar{\rho})$ to $\mathbb{Q}(\bar{\rho}_{k,\varepsilon})$. However this contradicts our choice of $B$ in the pair $(A, B)$ associated to $l_{k+1}$ through $\rho_k$. Thus there exists $f_{k+1} \in H^1(G_{S_k}, \mathrm{Ad}^0 \bar{\rho})$ such that $f_{k+1} \mid_{G_{l_k}}$ is unramified at $l_k$ and nontrivial. As the null 1-cohomology class at $l_{k+1}$ is ramified at $l_{k+1}$ the result follows.

LEMMA 6. *For $1 \leq i < k+1$ we have that $f_i \mid_{G_{l_{k+1}}}$ is trivial.*

*Proof.* The primes above $l_{k+1}$ in $\mathbb{Q}(\bar{\rho})$ split completely from $\mathbb{Q}(\bar{\rho})$ to $\mathbb{Q}(\bar{\rho}_{k-1,\varepsilon})$. Indeed, this is equivalent to the statement that the matrix $B$ has trivial projection from $N_k$ to $N_{k-1}$. Thus we see for any $g \in H^1(G_{S_{k-1}}, \mathrm{Ad}^0 \bar{\rho})$ that $g \mid_{G_{l_{k+1}}}$ is trivial. As $i < k+1$ we have $f_i \in H^1(G_{S_{i-1}}, \mathrm{Ad}^0 \bar{\rho}) \subset H^1(G_{S_{k-1}}, \mathrm{Ad}^0 \bar{\rho})$ and we are done.

*Remark.* In this inductive procedure, once chosen the $f_i$ remain fixed. For fixed $i$ the $u_{i,n}$ vary but $\lim_{n \to \infty} u_{i,n}$ exists in $\mathbb{Z}_p$.

PROPOSITION 2. *Let $\bar\rho$ be the mod $p$ representations coming from the elliptic curve $E$ as described in the introduction. Suppose $\{\rho_0, \rho_1, \ldots, \rho_k\}$ are surjective deformations of $\bar\rho$ to $\mathrm{GL}_2(\mathbb{Z}_p)$ such that for $1 \leq m \leq n \leq k$ that $\rho_n \mid_{G_{l_m}}$ is given by $\sigma_{l_m} \mapsto \begin{pmatrix} l_m & 0 \\ 0 & 1 \end{pmatrix}$ and $\tau_{l_m} \mapsto \begin{pmatrix} 1 & p^m u_{m,n} \\ 0 & 1 \end{pmatrix}$. Suppose also $\rho_{n-1} \equiv \rho_n \bmod p^n$ and there exists $f_n \in H^1(G_{S_{n-1}}, \mathrm{Ad}^0 \bar\rho)$ independent of the null 1-cohomology class at $l_n$ and $m \leq n \leq k+1$ implies $f_m \mid_{G_{l_n}}$ is trivial. Then there exists a prime $l_{k+1}$ and a surjective representation $\rho_{k+1} : G_{S_{k+1}} \to \mathrm{GL}_2(\mathbb{Z}_p)$ infinitely ramified at $l_{k+1}$ with $\rho_k \equiv \rho_{k+1} \bmod p^{k+1}$ and $\rho_{k+1} \mid_{G_{l_n}}$ is given by $\sigma_{l_n} \mapsto \begin{pmatrix} l_n & 0 \\ 0 & 1 \end{pmatrix}$ and $\tau_{l_n} \mapsto \begin{pmatrix} 1 & p^n u_{n,k+1} \\ 0 & 1 \end{pmatrix}$ for $1 \leq n \leq k+1$. Also, there exists $f_{k+1} \in H^1(G_{S_n}, \mathrm{Ad}^0 \bar\rho)$ independent of the null 1-cohomology class at $l_{k+1}$.*

*Proof.* Consider $\rho_k \bmod p^{k+2}$. This is ramified at $l_1, l_2, \ldots, l_k$. Introduce ramification at $l_{k+1}$ by adjusting $\rho_k \bmod p^{k+2}$ by *any* element of $H^1(G_{S_{k+1}}, \mathrm{Ad}^0 \bar\rho)$ ramified at $l_{k+1}$. This takes care of our first task.

Now we need to make sure that for each $i$, $k+1 \geq i \geq 1$ the local at $l_i$ deformation problem is unobstructed. We do this by forcing them to be special. First we do this for $l_{k+1}$ by adjusting by a suitable multiple of $f_{k+1}$ provided by Lemma 5. Then we adjust by a suitable multiple of $f_k$ for $l_k$, $f_{k-1}$ for $l_{k-1}$ and so on. As remarked above, since for $i > j$ we have $f_j \mid_{G_i}$ is trivial, adjusting by $f_j$ does not affect the deformation problem for $l_i$. We see that the lifting problem is locally unobstructed at all $l_i$ and thus globally unobstructed. We lift to mod $p^{n+2}$. Repeat this last process and lift to mod $p^{n+3}$. Continuing, we get our $\rho_{k+1}$ special at all primes of $\{l_1, \ldots, l_{k+1}\}$.

We have proved the following theorem.

THEOREM 1. *Fix $E_{/\mathbb{Q}}$ a semistable elliptic curve with good reduction at 3. For primes $p \geq 5$ in a set of density one there exist surjective representations $G_\mathbb{Q} \to \mathrm{GL}_2(\mathbb{Z}_p)$ ramified at infinitely many primes. The reduction mod $p$ of these representations is the Galois action on the $p$-torsion of $E$.*

*Proof.* Let $\rho$ be the limit of the $\rho_k$.

We now turn to our GRH results. We need a few preliminaries.

*Local at $S_0$ theories.* Let $v \in S_0$, $v \neq p$. Recall $E$ is semistable and thus has multiplicative reduction at $v$. Then from the theory of the Tate curve $\bar\rho \mid_{G_v} = \begin{pmatrix} \psi\chi & * \\ 0 & \psi^{-1} \end{pmatrix}$. The $*$ here is nontrivial as $p$ does not divide the $v$-adic valuation of $j(E)$.

LEMMA 7. *For all $v \in S_0$, $v \neq p$, $H^i(G_v, \mathrm{Ad}^0\bar{\rho})$ is trivial for $i = 0, 1, 2$.*

*Proof.* The $H^0$ result follows immediately. The $H^2$ result follows from local duality and requires $v \neq 3$. This is another reason why we insist that 3 be a prime of good reduction of $E$. The $H^1$ result follows from applying the local Euler characteristic, keeping in mind that $v$ is prime to $p$.

The deformation theory of $\bar{\rho}\mid_{G_v}$ is then trivial; i.e. the universal deformation ring for the local at $v$ problem is just $\mathbb{Z}_p$. For our purposes this means that there are no local at $v$ conditions in the global (ordinary at $p$) weight-2 Selmer group. See [W] for a discussion of Selmer groups.

We have chosen $p$ so $\bar{\rho}\mid_{G_p}$ is ordinary; that is

$$\bar{\rho}\mid_{G_p}(\alpha) = \begin{pmatrix} \psi\chi(\alpha) & * \\ 0 & \psi^{-1}(\alpha) \end{pmatrix}$$

where $\psi$ is an unramified character of order greater than 2 and $\chi$ is the cyclotomic character. The first requirement corresponds to the fact that $a_p \neq \pm 1$ for our elliptic curve $E/\mathbb{Q}$ and guarantees that $H^2(G_p, \mathrm{Ad}^0\bar{\rho}) = 0$ (see [M2]). Under these circumstances Wiles and Taylor-Wiles have proved that the minimal global universal ordinary at $p$ weight-2 deformation ring is just $\mathbb{Z}_p$. See [Da] for an explicit example involving the elliptic curve $X_0(17)$ and the prime $p = 5$.

LEMMA 8. *$H^2(G_p, \mathrm{Ad}^0\bar{\rho}) = 0$. $H^0(G_p, \mathrm{Ad}^0\bar{\rho})$ is zero- or one-dimensional as the $*$ in $\bar{\rho}\mid_{G_p}$ is nontrivial or trivial. $H^1(G_p, \mathrm{Ad}^0\bar{\rho})$ is three- or four-dimensional as the $*$ is nontrivial or trivial.*

*Proof.* The result for $H^2$ follows from local duality using that $\psi^2$ is not the trivial character. The $H^0$ result is immediate and the $H^1$ result is a consequence of the local Euler-Poincaré characteristic.

LEMMA 9. *$H^1_{\mathrm{ord}}(G_p, \mathrm{Ad}^0\bar{\rho})$ is one- or two-dimensional as the $*$ is nontrivial or trivial.*

*Proof.* $H^1_{\mathrm{ord}}$ consists of the 1-cohomology classes that give rise to ordinary deformations to the dual numbers. If $*$ is nontrivial, then in Section 6 of [R3] all reducible lifts of $\bar{\rho}$ to $\mathbb{Z}/p^2$ are computed. There are $p$ of these that are ordinary and have determinant the cyclotomic character so in this case $H^1_{\mathrm{ord}}$ is one-dimensional.

If the $*$ is trivial $H^1(G_p, \mathrm{Ad}^0\bar{\rho})$ is four-dimensional by Lemma 8. In this case we can write down the ordinary deformation to the dual numbers $\tilde{\rho}$. One is given by $\tilde{\rho}\mid_{G_p}(\alpha) = \begin{pmatrix} \psi\chi(\alpha) & \varepsilon* \\ 0 & \psi^{-1}(\alpha) \end{pmatrix}$. There are also the ordinary

deformations given by

$$\tilde{\rho}|_{G_p}(\alpha) = \begin{pmatrix} \psi\chi(1+\varepsilon h)(\alpha) & 0 \\ 0 & \psi^{-1}(1-\varepsilon h)(\alpha) \end{pmatrix}$$

where $h$ is an unramified character of order $p$. It is straightforward to see these generate the ordinary local tangent space.

LEMMA 10. *The restriction map $H^1(G_{S_0}, \mathrm{Ad}^0\bar{\rho}) \to H^1(G_p, \mathrm{Ad}^0\bar{\rho})$ is injective. The image of this map is a two-dimensional space whose intersection with $H^1_{\mathrm{ord}}(G_p, \mathrm{Ad}^0\bar{\rho})$, the ordinary 1-cohomology classes, is trivial.*

*Proof.* Since we assume $H^2(G_{S_0}, \mathrm{Ad}^0\bar{\rho}) = 0$, global duality and the fact that $\bar{\rho}$ is odd imply $H^1(G_{S_0}, \mathrm{Ad}^0\bar{\rho})$ is two-dimensional. Recall, as mentioned before Lemma 8, that the universal ordinary at $p$ weight-2 deformation ring is just $\mathbb{Z}_p$. If the restriction map $H^1(G_{S_0}, \mathrm{Ad}^0\bar{\rho}) \to H^1(G_p, \mathrm{Ad}^0\bar{\rho})$ were not injective then the kernel of this map would give rise to (trivial) ordinary at $p$ lifts to the dual numbers of $\mathbb{F}_p$; that is the universal ordinary at $p$ weight-2 ring would be a nontrivial quotient of $\mathbb{Z}_p[[T_1, T_2, \ldots, T_r]]$ for some $r > 0$, not just $\mathbb{Z}_p$. Thus the restriction map is injective. That the image intersects $H^1_{\mathrm{ord}}(G_p, \mathrm{Ad}^0\bar{\rho})$ trivially follows similarly.

COROLLARY 1. $H^1(G_p, \mathrm{Ad}^0\bar{\rho}) \simeq H^1_{\mathrm{ord}}(G_p, \mathrm{Ad}^0\bar{\rho}) \oplus \mathrm{res}_{G_p}(H^1(G_{S_0}, \mathrm{Ad}^0\bar{\rho}))$.

*Proof.* Since the image of the restriction map intersects the local ordinary tangent space trivially, counting dimensions in Lemmas 8, 9, and 10 proves the corollary.

COROLLARY 2. *Let $S_k$ be a set of primes containing $S_0$. Suppose $\psi_n$ is a global deformation of $\bar{\rho}$ to $\mathbb{Z}/p^n$ ramified only at primes in $S_k$ that is ordinary at $p$. Suppose there exists a deformation $\psi_{n+1}$ of $\psi_n$ to $\mathbb{Z}/p^{n+1}$. There exists $f \in H^1(G_{S_0}, \mathrm{Ad}^0\bar{\rho})$ such that $f.\psi_{n+1}$ is ordinary at $p$.*

*Proof.* Clearly there is an ordinary deformation of the local at $p$ representation $\psi_n|_{G_p}$ to $\mathrm{GL}_2(\mathbb{Z}/p^{n+1})$. Thus there exists a $g \in H^1(G_p, \mathrm{Ad}^0\bar{\rho})$ such that $g.\psi_{n+1}|_{G_p}$ is ordinary. By Corollary 1 we can uniquely write $g = f + h$ where $f$ is in the image of the the restriction map of $H^1(G_S, \mathrm{Ad}^0\bar{\rho})$ in $H^1(G_p, \mathrm{Ad}^0\bar{\rho})$ and $h \in H^1_{\mathrm{ord}}(G_p, \mathrm{Ad}^0\bar{\rho})$. Thus there exists a global 1-cohomology class $\tilde{f} \in H^1(G_S, \mathrm{Ad}^0\bar{\rho})$ with $\mathrm{res}_p(\tilde{f}) = f$. We see $\tilde{f}.\psi_{n+1}|_{G_p} = f.\psi_{n+1}|_{G_p} = (g-h).\psi_{n+1}|_{G_p}$. Since $h \in H^1_{\mathrm{ord}}(G_p, \mathrm{Ad}^0\bar{\rho})$ and $g.\psi_{n+1}|_{G_p}$ is ordinary, so is $(g-h).\psi_{n+1}|_{G_p}$.

For $a$ nonzero in $\mathbb{Z}/p$ let $a^*$ denote the Teichmüller lift of $a$ to $\mathbb{Z}_p$, i.e. the unique $p-1^{\mathrm{st}}$ root of unity in $\mathbb{Z}_p$ congruent to $a$ mod $p$. As $p \geq 5$ we see $2^* \neq (1/2)^*$. Let $A \in \mathrm{GL}_2(\mathbb{Z}/p^{k+2})$ be the matrix $\begin{pmatrix} 2^*(1+p^{k+1}) & 0 \\ 0 & 1-p^{k+1} \end{pmatrix}$.

Let $\mathbf{0} = \begin{pmatrix} 0 & 0 \\ 0 & 0 \end{pmatrix}$ be the trivial element of $N_k$. We give the trivial element a matrix description to emphasize the $\mathrm{Gal}(\mathbb{Q}(\bar{\rho})/\mathbb{Q})$ action.

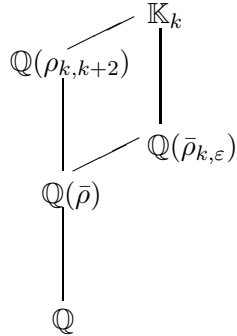

Recall that we denoted $\mathrm{GL}_2(\mathbb{Z}/p^{k+2}) \simeq \mathrm{Gal}(\mathbb{Q}(\bar{\rho}_{k,k+2})/\mathbb{Q})$ by $C_k$ and $\mathrm{Gal}(\mathbb{Q}(\bar{\rho}_{k,\varepsilon})/\mathbb{Q}(\bar{\rho})) \simeq \mathrm{Gal}(\mathbb{K}_k/\mathbb{Q}(\bar{\rho}_{k,k+2}))$ by $N_k$ and $\mathrm{Gal}(\mathbb{K}_k/\mathbb{Q})$ by $H_k$. Lemma 3 showed $H_k$ was a semidirect product of $N_k$ by $C_k$.

The lemma below follows immediately from Chebotarev's theorem.

LEMMA 11. *If $p \geq 5$ then there are infinitely many primes whose Frobenius in $\mathrm{Gal}(\mathbb{K}_k/\mathbb{Q})$ is in the conjugacy class of $(A, \mathbf{0})$ in the semi direct product of $C_k$ by $N_k$.*

Note that for such primes $q$ we have $q \equiv \det(\rho_k(A)) \equiv 2^* \bmod p^{k+2}$. Our choice for $A$ instead of $\begin{pmatrix} 2^* & 0 \\ 0 & 1 \end{pmatrix}$ is so that we can guarantee that we add ramification at our new prime(s) exactly mod $p^{k+2}$. Also note that unlike the unconditional situation, for such primes $q$, by our choice of the matrix $\mathbf{0}$, the primes of $\mathbb{Q}(\bar{\rho})$ above $q$ split completely from $\mathbb{Q}(\bar{\rho})$ to $\mathbb{Q}(\bar{\rho}_{k,\varepsilon})$.

*Increasing the ramification.* We induct as before. For technical reasons we may have to add two primes at a time to the level. The prime(s) at which we wish to allow ramification will be as in Lemma 11.

Suppose now for $0 \leq n \leq k$ that $S_n = S_{n-1} \cup X_n$ where $X_n = \{l_n\}$ or $\{l_{n1}, l_{n2}\}$. Suppose also that for such $n$ we have constructed $\rho_n : G_{S_n} \to \mathrm{GL}_2(\mathbb{Z}_p)$ surjective and ramified at all primes in $S_n$. Also, we assume $\rho_{n-1} \equiv \rho_n$ mod $p^n$. Suppose further for each prime $q \in X_n$ there exists a global 1-cohomology class $f_q \in H^1(G_{S_n}, \mathrm{Ad}^0 \bar\rho)$ that is *independent of the null 1-cohomology class at $q$*. Furthermore, for $q \in X_n$, $r \in X_m$ and $n < m$ we assume that $f_q |_{G_r}$ is trivial. We will complete the induction by constructing $X_{k+1}, \rho_{k+1}$ and for each $q \in X_{k+1}$ the global 1-cohomology class $f_q$.

We have the exact sequence

$$0 \to \mathrm{III}^2_{S_{k+1}}(\mathrm{Ad}^0 \bar\rho) \to H^2(G_{S_{k+1}}, \mathrm{Ad}^0 \bar\rho) \to \oplus_{v \in S_{k+1}} H^2(G_v, \mathrm{Ad}^0 \bar\rho).$$

As in the argument following Lemma 4, $\mathrm{III}^2_{S_{k+1}}(\mathrm{Ad}^0 \bar\rho)$ is trivial. Thus we need only analyze local deformation problems to analyze global deformation problems. Since for $v \in S_0$ we are assuming that $H^2(G_v, \mathrm{Ad}^0 \bar\rho)$ is trivial we only study primes in $S_{k+1} - S_0$. Indeed, the right-hand map above is an isomorphism.

Consider the primes $q$ in our Chebotarev class of Lemma 11. These $q$ are candidates for the prime(s) we would like to add to the level. Let us consider elements

$$f \in H^1(G_{S_k \cup \{q\}}, \mathrm{Ad}^0 \bar\rho) - H^1(G_{S_k}, \mathrm{Ad}^0 \bar\rho).$$

Recall from Lemma 1 that $H^2(G_q, \mathrm{Ad}^0 \bar\rho)$ is one-dimensional. By Fact 1 we see $H^1(G_{S_k}, \mathrm{Ad}^0 \bar\rho)$ is of codimension 1 in $H^1(G_{S_k \cup \{q\}}, \mathrm{Ad}^0 \bar\rho)$. We need for $f |_{G_q}$ to be independent of the null 1-cohomology class for the prime $q$. Note that $f$ is clearly ramified at $q$. The choice of **0** in Lemma 11 guarantees the primes above $q$ in $\mathbb{Q}(\bar\rho)$ split completely from $\mathbb{Q}(\bar\rho)$ to $\mathbb{Q}(\bar\rho_{k,\varepsilon})$. Thus for any $g \in H^1(G_{S_k}, \mathrm{Ad}^0 \bar\rho)$ we have that $g |_{G_q}$ is trivial. The nonnullity of $f |_{G_q}$ is not affected by adding to $f$ some $g \in H^1(G_{S_k}, \mathrm{Ad}^0 \bar\rho)$ and is therefore a well-defined notion. If for some prime $q$ in our Chebotarev class we have $f |_{G_q}$ is nonnull then we can choose $l_{k+1}$ to be this $q$, $X_{k+1} = \{l_{k+1}\}$ and $f_{k+1}$ to be this $f$. We start at mod $p^{k+2}$, that is we consider $\rho_k$ mod $p^{k+2}$. Note $f_{k+1} |_{G_{l_{k+1}}}$ is ramified and that $\rho_k(\mathrm{Frob}_{l_{k+1}}) \equiv \begin{pmatrix} l_k(1+p^{k+1}) & 0 \\ 0 & 1 - p^{k+1} \end{pmatrix}$ mod $p^{k+2}$ is *not* special. (But it is special mod $p^{k+1}$.) Thus at the mod $p^{k+2}$ stage we *must* alter it by a nonzero multiple of $f_{k+1}$ to get to the special form. This introduces ramification at $l_{k+1}$ exactly mod $p^{k+2}$.

We now alter the deformation problems at prime(s) of $X_k, X_{k-1}, ... X_2, X_1$ successively by multiples of $f_{l_i}$ (or by $f_{l_{i1}}$ and $f_{l_{i2}}$ if $X_i$ has two elements) for $l_i \in X_i$ to remove obstructions at these problems. The only difficulty is that for $i > j$ altering by $f_j$ may introduce an obstruction to the deformation problem at $r \in X_i$. This does not happen because by construction $i > j$ implies for $w \in X_j, r \in X_i$ that $f_w |_{G_r}$ is trivial. Finally, by Corollary 2 we can use an element of $H^1(G_{S_0}, \mathrm{Ad}^0 \bar\rho)$ to force the local at $p$ problem to be ordinary. Again, this will not change any of the local at $l_i$ obstructions. We can lift to mod $p^{k+3}$, knowing we have ramification at $l_{k+1}$. This deformation may be obstructed at any of the $l_k$ and nonordinary at $p$. To fix these problems successively adjust by appropriate multiples of $f_{l_{k+1}}, f_{l_k}, ... f_{l_2}, f_{l_1}$ (or by $f_{l_{i_1}}$ and $f_{l_{i_2}}$ if $X_i$ has two elements) to remove obstructions at the $l_k$. Then adjust by a suitable element of $H^1(G_{S_0}, \mathrm{Ad}^0 \bar{\rho})$ to guarantee the deformation is ordinary at $p$. Now lift to mod $p^{k+4}$. We continue this process to characteristic zero and use Wiles and Taylor-Wiles to get $\rho_{k+1}$ modular of level prime to $p$.

In the unconditional construction we used an element of $H^1(G_{S_0}, \mathrm{Ad}^0 \bar{\rho})$ as our nonnull 1-cohomology class for the prime $l_1$. Thus we had only one dimension of $H^1(G_{S_0}, \mathrm{Ad}^0 \bar{\rho})$ "left over". By the above arguments, particularly the application of Corollary 2, we see this is not enough to arrange that our representations be ordinary at $p$.

We now turn to the case where we assume all primes $q$ with Frobenius in the Chebotarev class of Lemma 11 are null; that is for

$$f_q \in H^1(G_{S_k \cup \{q\}}, \mathrm{Ad}^0 \bar{\rho}) - H^1(G_{S_k}, \mathrm{Ad}^0 \bar{\rho})$$

we assume $f_q \mid_{G_q}$ is null.

*Analytic considerations.* Now suppose that for *all* $q$ in our Chebotarev class that the corresponding 1-cohomology class

$$f_q \in H^1(G_{S_k \cup \{q\}}, \mathrm{Ad}^0 \bar{\rho}) - H^1(G_{S_k}, \mathrm{Ad}^0 \bar{\rho})$$

is null at $q$. (This seems extremely unlikely! The most naive heuristic suggests that these 1-cohomology classes are null with probability $1/p$. The little data that has been gathered seems to confirm this.) We get around this problem by adding *two* primes to the level simultaneously, but we need the GRH.

For each $q_i$ in the Chebotarev class let $\mathbb{K}_k^{q_i}$ denote the composite of $\mathbb{K}_k$ and $\mathbb{Q}(\bar{\rho}_{k,q_i,\varepsilon})$, the field cut out by lifts to the dual numbers of $\bar{\rho}$ factoring through $G_{S_k \cup \{q_i\}}$. We now ask how $q_j$ splits in $\mathbb{Q}(\bar{\rho}_{k,q_i,\varepsilon})$, i.e. how primes of $\mathbb{Q}(\bar{\rho})$ above $q_j$ split from $\mathbb{Q}(\bar{\rho})$ to $\mathbb{Q}(\bar{\rho}_{k,q_i,\varepsilon})$. Recall we are assuming $f_{q_i} \mid_{G_{q_i}}$ is null.

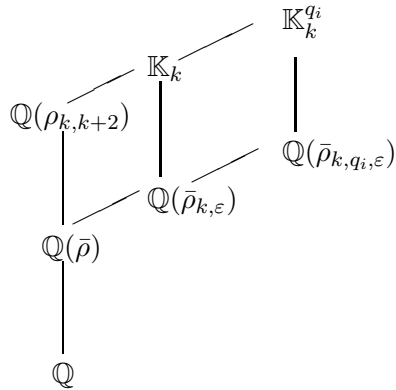

Recall $H_k = \mathrm{Gal}(\mathbb{K}_k/\mathbb{Q})$, $N_k = \mathrm{Gal}(\mathbb{Q}(\bar{\rho}_{k,\varepsilon})/\mathbb{Q}(\bar{\rho}))$, and

$$C_k = \mathrm{Gal}(\mathbb{Q}(\rho_{k,k+2})/\mathbb{Q}) \simeq \mathrm{GL}_2(\mathbb{Z}/p^{k+2}).$$

By Lemma 3, $H_k$ is the semidirect product of $C_k$ by $N_k$.

Let
$$B_k^{q_i} = \mathrm{Gal}(\mathbb{K}_k^{q_i}/\mathbb{K}_k) \simeq \mathrm{Gal}(\mathbb{Q}(\bar{\rho}_{k,q_i,\varepsilon})/\mathbb{Q}(\bar{\rho}_{k,\varepsilon}))$$

and $D_k^{q_i} = \mathrm{Gal}(\mathbb{K}_k^{q_i}/\mathbb{Q})$. Then $D_k^{q_i}$ is the semidirect product of $B_k^{q_i}$ by $H_k$. The proof is the same as in Lemma 3. Note that $H_k$ acts on $B_k^{q_i}$ and this action factors through the quotient $\mathrm{Gal}(\mathbb{Q}(\bar{\rho})/\mathbb{Q})$. Furthermore, essentially by Fact 1, $B_k^{q_i} \simeq \mathrm{Ad}^0 \bar{\rho}$ where the isomorphism is as $\mathbb{F}_p[H_k]$-modules. Note the $H_k$ action on $B_k^{q_i}$ factors through a quotient independent of $k$. The isomorphism classes of the *groups* $B_k^{q_i}$, and $D_k^{q_i}$ depend only on $k$ and *not* on $q_i$.

Consider the conjugacy class
$$C = \left( \begin{pmatrix} 2^*(1+p^{k+1}) & 0 \\ 0 & 1-p^{k+1} \end{pmatrix}, \mathbf{0} \right)$$

in $H_k$ the semidirect product of $C_k$ by $N_k$. Since $D_k^{q_i}$ is a semidirect product of $B_k^{q_i}$ by $H_k$ we may extend $C$ to a conjugacy class

$$\tilde{C}^{q_i} = (C, \mathbf{0}) = \left( \left( \begin{pmatrix} 2^*(1+p^{k+1}) & 0 \\ 0 & 1-p^{k+1} \end{pmatrix}, \mathbf{0} \right), \mathbf{0} \right)$$

in $D_k^{q_i}$, the semidirect product of $B_k^{q_i}$ by $H_k$.

LEMMA 12. *Let $d$ be the density of primes with Frobenius in $C$ in $H_k$. Then the density of primes with Frobenius in $\tilde{C}^{q_i}$ in $D_k^{q_i}$ is $d/p$ and the density of primes whose Frobenius in $D_k^{q_i}$ projects to an element of $C$ but is not in $\tilde{C}^{q_i}$ is $d(1-1/p)$.*

*Proof.* The density of primes with Frobenius in a given conjugacy class is the reciprocal of the order of the centralizer of an element in the class.

So we must compare the order of the centralizer of an element of $C$ in $H_k$ with that of an element of $\tilde{C}^{q_i}$ in $D_k^{q_i}$. But $D_k^{q_i}$ is "bigger" than $H_k$ by a copy of $\mathrm{Ad}^0 \bar{\rho}$ with action through $\bar{\rho}$. The centralizer of an element of $\tilde{C}^{q_i}$ in $\mathrm{Ad}^0 \bar{\rho}$ is just the set of trace 0 diagonal matrices. These have order $p$ so our density is $d/p$. The other statements are immediate.

Consider the infinite matrix below:

|  | $q_1$ | $q_2$ | ... | $q_r$ | ... |
|---|---|---|---|---|---|
| $\mathbb{K}_k^{q_1}$ | $*$ | 0 | 1 | 1 | 1... |
| $\mathbb{K}_k^{q_2}$ | 1 | $*$ | 1 | 0 | 1... |
| ... | 1 | 1 | $*$ | 0 | 1... |
| ... | .. | .. | .. | .. | .. |
| $\mathbb{K}_k^{q_r}$ | .. | .. | .. | .. | .. |
| ... | .. | .. | .. | .. | .. |

A "0" in the $ij$ spot indicates that $q_j$ has Frobenius $\tilde{C}^{q_i} = (C, \mathbf{0})$ in the semidirect product $D_k^{q_i}$ of $B_k^{q_i}$ by $H_k$. Recall $B_k^{q_i} = \mathrm{Gal}(\mathbb{Q}(\bar{\rho}_{k,q_i,\varepsilon})/\mathbb{Q}(\bar{\rho}_{k,\varepsilon}))$ and as an $\mathbb{F}_p[\mathrm{Gal}(\mathbb{Q}(\bar{\rho})/\mathbb{Q})]$-module is isomorphic to $\mathrm{Ad}^0 \bar{\rho}$. Thus a "0" in the $ij$ spot indicates that the primes above $q_j$ in $\mathbb{Q}(\bar{\rho})$ split completely in $\mathbb{Q}(\bar{\rho}_{k,q_i,\varepsilon})$. A "1" in the $ij$ spot means the primes above $q_j$ in $\mathbb{Q}(\bar{\rho})$ do not split completely in $\mathbb{Q}(\bar{\rho}_{k,q_i,\varepsilon})$.

LEMMA 13. *Suppose in the above matrix there are "1"'s in the $ij$ and $ji$ entries. Then $f_{q_i}\mid_{G_{q_j}}$ and $f_{q_j}\mid_{G_{q_i}}$ are unramified and nontrivial.*

*Proof.* Recall $f_{q_i} \in H^1(G_{S_k \cup \{q_i\}}, \mathrm{Ad}^0 \bar{\rho}) - H^1(G_{S_k}, \mathrm{Ad}^0 \bar{\rho})$. Clearly $f_{q_i}\mid_{G_{q_j}}$ is unramified. Since $q_j$ (and $q_i$) have Frobenius in the conjugacy class $C$ of the group $H_k$, the primes above $q_j$ in $\mathbb{Q}(\bar{\rho})$ split completely from $\mathbb{Q}(\bar{\rho})$ to $\mathbb{Q}(\bar{\rho}_{k,\varepsilon})$. If $f_{q_i}\mid_{G_{q_j}}$ were trivial the primes above $q_j$ would split completely from $\mathbb{Q}(\bar{\rho})$ to $\mathbb{Q}(\bar{\rho}_{k,q_i,\varepsilon})$. The "1" in the $ij$ entry prevents this. The same argument works for $ji$.

The condition of Lemma 13 is exactly what we need to add both $q_i$ and $q_j$ to the level simultaneously. If this condition is met then we let $l_{k+1,1} = q_i$ and $l_{k+1,2} = q_j$.

Suppose now that a "1" in the $ij$ spot always implies the $ji$ spot has a "0". We make a naive argument to show this is not reasonable and then prove this assuming the GRH.

Recall $d$ is the density of primes in the Chebotarev class $C$ of the group $H_k$. Let us ask how many "1"'s we expect in row $i$ for $q_j < x$. The "0"'s in row $i$ correspond to primes with Frobenius in $\tilde{C}^{q_i}$. Thus the "1"'s correspond to primes whose Frobenius in $H_k$ lies in the conjugacy class $C$, but whose Frobenius in $D_k^{q_i}$ does *not* lie in the class $\tilde{C}^{q_i}$. Lemma 12 and Chebotarev's theorem say this number is $d(1-1/p)li(x) + E_j(x)$ where $E_j(x)$ is an error term as are subsequent $E$ terms and $li(x)$ is the logarithmic integral, $\int_2^x dt/\log(t)$. Summing over all the rows up to $x$ (roughly $dli(x)$ primes) we see that the number of "1"'s in a big $x \times x$ square is $d^2(1 - 1/p)li^2(x) + E(x)$.

Let us find an upper bound for the number of "1"'s in our $x \times x$ square by finding an upper bound for the number of "1"'s in each column and summing over the columns. Since we are assuming that row$(i) + $ col$(i)$ has only "0"'s and "1"'s (i.e. no "2"'s) we may assume that col$(i) \leq \vec{1} - $ row $i$ where $\vec{1}$ denotes the row of all "1"'s. Thus the number of "1"'s in column $i$ is, by Chebotarev's theorem, less than $(d/p)li(x) + \tilde{E}_i(x)$. That we can use Chebotarev's theorem to estimate the number of "1"'s in a column from above follows from our assumption that a "1" in the $ij$ spot implies there is a "0" in the $ji$ spot. Summing over all the columns up to $x$ we see that the number of "1"'s in the $x \times x$ square is at most $(d^2/p)li^2(x) + \tilde{E}(x)$. Since this last sum is an upper bound for the number of "1"'s in an $x \times x$ square we have $(d^2/p)li^2(x) + \tilde{E}(x) \geq d^2(1-1/p)li^2(x) + E(x)$. We see that $\tilde{E}(x) - E(x) > d^2(1 - 2/p)li^2(x)$. At a minimum this seems unlikely. We show using estimates of Lagarias and Odlyzko (which assume the GRH) that the error term $\tilde{E}(x) - E(x)$ is $o(x^{3/2+\varepsilon})$. Thus choosing $x$ suitably large will give a contradiction.

LEMMA 14. *The absolute discriminant of $\mathbb{K}_k^{q_i}$ is $c_1 q_i^{c_2}$ where $c_1$ and $c_2$ are constants depending only on $k$.*

*Proof.* As $\text{Disc}(\mathbb{K}_k/\mathbb{Q})$ is independent of $q_i$ we need only study $\text{Disc}(\mathbb{K}_k^{q_i}/\mathbb{K}_k)$. Its contribution to $\text{Disc}(\mathbb{K}_k^{q_i}/\mathbb{Q})$ will be a power of $\text{Disc}(\mathbb{K}_k^{q_i}/\mathbb{K}_k)$ that depends only on $k$. Note that $\mathbb{K}_k^{q_i}/\mathbb{K}_k$ is only ramified at primes in $S_k \cup \{q_i\}$ and $[\mathbb{K}_k^{q_i} : \mathbb{K}_k] = p^3$. At all primes of $\mathbb{K}_k$ except those above $p$ the ramification is tame and the contribution to $\text{Disc}(\mathbb{K}_k^{q_i}/\mathbb{K}_k)$ is at most a (universally) bounded power in these primes. At $p$ the extension corresponds to a lift to the dual numbers so we may compute the (local at $p$) discriminant of *all* such lifts and treat this as a constant. As we have bounded $\text{Disc}(\mathbb{K}_k^{q_i}/\mathbb{K}_k)$ appropriately the result follows.

We recall the theorem of [LO].

THEOREM B. *Assume the GRH. Let $\mathbb{L}/\mathbb{K}$ be a Galois extension of number fields with Galois group $H$. Let $C$ be a conjugacy class of $H$. Denote by $\pi_C(X)$ the number of prime ideals of $\mathbb{K}$ of norm less than $x$ whose Frobenii lie in $C$. Then $\left|\pi_C(x) - \frac{|C|}{|H|}li(x)\right| \leq e_1\left(\frac{|C|}{|H|}x^{1/2}\log(D_\mathbb{L} x^{n_\mathbb{L}}) + \log(D_\mathbb{L})\right)$ for all $x > 2$.*

Here $D_\mathbb{L}$ is the absolute discriminant of $\mathbb{L}$ and $n_\mathbb{L}$ is the degree $[\mathbb{L} : \mathbb{Q}]$ and $e_1$ is an absolute constant.

Note that Lagarias and Odlyzko have an unconditional result. If one only considers cases where there is no Seigel zero, which should suffice for some of our applications, then the error term in the unconditional result is manageable. The difficulty is that their estimate only holds for $x$ much bigger than $D_\mathbb{L}$. As we want to consider $q_i$ up to $x$ and let $\mathbb{K}_k^{q_i}$ play the role of $\mathbb{L}$ we cannot apply their unconditional result. Put $n := [\mathbb{K}_k^{q_i} : \mathbb{Q}]$. Note $n$, $c_1$ and $c_2$ are constant as the $q_i$ vary through our Chebotarev class.

PROPOSITION 3. *Assume the GRH. In the matrix above there exist (infinitely many) pairs of integers $(i, j)$ such that there are "1"'s in both the $ij$ and $ji$ spots.*

*Proof.* Assume the proposition is false. We sum the "1"'s in an $x \times x$ square as described above. Recall the density of the primes $\{q_i\}$ in $C$ is $d$. Then by Chebotarev's theorem in row $i$ the density of "1"'s is $d(1 - 1/p)$. The theorem of Lagarias and Odlyzko above shows that we expect at most

$$d(1 - 1/p)li(x) + e_1\left(d(1 - 1/p)x^{1/2}\log(c_1 q_i^{c_2} x^n) + \log(c_1 q_i^{c_2})\right)$$

"1"'s, where $c_1$ and $c_2$ are as in Lemma 14 and $n = [\mathbb{K}_k^{q_i} : \mathbb{Q}]$. Keeping in mind that the density of the $\{q_i\}$ is $d$, summing over all $q_i < x$ gives

$$\sum_{q_i < x, q_i \in C} \left(d(1 - 1/p)li(x)\right) + O\left(e_1 d(1 - 1/p)x^{1/2} n \log(x)\right)$$
$$+ O\left(e_1 d(1 - 1/p)x^{1/2} \log(c_1)\right) + O\left(e_1 d(1 - 1/p)x^{1/2} c_2 \log(q_i)\right)$$
$$+ O\left(e_1 \log(c_1)\right) + O\left(e_1 c_2 \log(q_i)\right)$$

"1"'s in our $x \times x$ square. The main term for the number of "1"'s is $d^2(1 - 1/p)li^2(x)$. Recall the well-known sum $\sum_{q_i \leq x, q_i \in C} \log(q_i) = O(x)$. Then the error term in our sum is

$$O(li(x)) \cdot O(x^{1/2} \log x) + O(x^{1/2} \log(x)) \cdot O(li(x)) + O(x^{1/2}) \cdot O(li(x))$$
$$+ O(x^{1/2}) \cdot O(x) + O(1) \cdot O(li(x)) + O(1) \cdot O(x) = o(x^{3/2 + \varepsilon}).$$

Note that the constants in the $O$ terms depend on $n = [\mathbb{K}_k^{q_i} : \mathbb{Q}]$ which in turn depends on $k$, and on the density $d$ which also depends on $k$. Thus the constants depend only on $k$ which is fixed.

Let us now sum the "1"'s by columns, assuming that column $i$ is dominated by row $i$. We estimate this sum from above, assuming that column $i = \vec{1}-$ row $i$. The number of "1"'s in column $i$ is then less than $d\pi_C(x) -$ (the number of "1"'s in row $i$) which we have computed above. Summing over $q_i \leq x$, we get the number of "1"'s in our square to be $(d^2/p)li^2(x) +$ an error term the same size as before. Equating, we see $d^2(1 - 2/p)li^2(x) = o(x^{3/2+\varepsilon})$ which is false for large $x$. If there were only finitely many pairs they could not "make up" for the discrepancies between the two summing methods. The proposition is proved.

PROPOSITION 4. *Assume that for all $q$ in our Chebotarev class the global 1-cohomology class $f \in H^1(G_{S_k \cup \{q\}}, \mathrm{Ad}^0 \bar{\rho}) - H^1(G_{s_k}, \mathrm{Ad}^0 \bar{\rho})$ is null at $q$. Then, assuming the* GRH, *there exist $q_i, q_j$ in our class at which we can allow ramification simultaneously to complete our induction.*

*Proof.* Use $q_i, q_j$ as given by Proposition 3. Let

$$h_1 \in H^1(G_{S_k \cup \{q_i\}}, \mathrm{Ad}^0 \bar{\rho}) - H^1(G_{S_k}, \mathrm{Ad}^0 \bar{\rho}),$$
$$h_2 \in H^1(G_{S_k \cup \{q_j\}}, \mathrm{Ad}^0 \bar{\rho}) - H^1(G_{S_k}, \mathrm{Ad}^0 \bar{\rho}).$$

Note that $h_1 \mid_{G_{q_i}}$ is ramified and null at $q_i$ and $h_2 \mid_{G_{q_j}}$ is ramified and null at $q_j$.

We will let $l_{k+1,1} = q_i$ and $l_{k+1,2} = q_j$, let $f_{k+1,1} = h_2$, $f_{k+1,2} = h_1$ and $X_{k+1} = \{l_{k+1,1}, l_{k+1,2}\}$. We claim that $f_{k+1,i}$ is nonnull for $l_{k+1,i}$ for $i = 1, 2$. Indeed, this follows from Proposition 3 and Lemma 13. Now consider $\rho_k \mod p^{k+2}$. First adjust the deformation problem at $l_{k+1,1}$ by a suitable multiple of $f_{k+1,1}$ and then adjust the deformation problem at $l_{k+1,2}$ by a suitable multiple of $f_{k+1,2}$. The only difficulty is that this last adjustment may cause an obstruction at $l_{k+1,1}$. But $f_{k+1,2} = h_1$ is null at $l_{k+1,1}$ and so causes no such problems.

Since $\rho_k(\text{Frob}_{l_{k+1,i}})$ has mod $p^{k+2}$ eigenvalues that are *not* $\{l_{k+1,i}, 1\}$ the multiples of $f_{k+1,i}$ that we adjust by are nonzero for $i = 1, 2$. This guarantees that we introduce ramification at these two new primes exactly at mod $p^{k+2}$.

Now we proceed as before, fixing the local deformation problems at primes of $X_k, X_{k-1}, ... X_2, X_1$ and the ordinariness problem at $p$. Lift to mod $p^{k+3}$ and repeat this procedure. We have proved Theorem 2 below.

THEOREM 2. *Assume the* GRH. *Fix $E_{/\mathbb{Q}}$ a semistable elliptic curve. For primes $p \geq 5$ in a set of density one and for every nonnegative integer $k$ there exist surjective representations $\rho_k : G_{S_k} \to \text{GL}_2(\mathbb{Z}_p)$ ramified at all primes of $S_k = S_0 \cup \cup_{i=1}^k X_k$ where $X_k$ is a set containing 1 or 2 primes. Each $\rho_k$ is modular of weight 2 and level prime to $p$. Furthermore, $\rho_{k-1} \equiv \rho_k \mod p^k$ and $\rho = \lim \rho_k$ exists and is ramified at all primes in all $S_k$ and is crystalline at $p$. The reduction $\mod p$ of these representations is the Galois action on the $p$-torsion of $E$.*

Our $\rho_k$ correspond to newforms $h_k$ of weight 2, trivial character, and level prime to $p$. The $q$-expansions of the $h_k$ converge $p$-adically. Diamond has observed that only finitely many of the $\rho_k$ come from elliptic curves over $\mathbb{Q}$. Indeed, if infinitely many came from elliptic curves over $\mathbb{Q}$ then for large $N$ we would have infinitely many elliptic curves with (Galois) isomorphic $p^N$ torsion. Thus we would have infinitely many points on a twist of the modular curve $X(p^N)$ rational over some number field. As $X(p^N)$ has genus bigger than 1 for large $N$ this is impossible by Faltings' theorem. Thus our representations correspond to weight-2 forms such that $p$ splits completely in the field generated by the all the eigenvalues of Frobenii.

*Remark.* Why did we choose $A = \begin{pmatrix} 2^*(1+p^{k+1}) & 0 \\ 0 & 1-p^{k+1} \end{pmatrix}$ as opposed to $\begin{pmatrix} 2^* & 0 \\ 0 & 1 \end{pmatrix}$? Either choice can be made to work. The point is that we want to guarantee that at mod $p^r$ for some $r$ we will actually adjust by a *nonzero* multiple of this global 1-cohomology class and introduce ramification

4814 RAVI RAMAKRISHNAat the new prime(s). The possibility that $\rho_{n-1}(\text{Frob}_{l_n})$ has eigenvalues $l_n$ and 1 means that we are never able to introduce ramification at $l_n$. (Of course in the odd modular case we do not expect this!) This is similar to Proposition 1 of [R2]. The first choice for $A$ above guarantees that we could introduce this ramification for $l_{k+1}$ exactly at mod $p^{k+2}$.

We now naively ask whether we expect to find one global 1-cohomology class for the two tasks at hand for the prime $l_{k+1}$ if we choose $A$ as in the GRH section of this paper. We need this ramified at $l_{k+1}$ global 1-cohomology class to have restriction at $l_{k+1}$ independent of the local at $l_{k+1}$ null 1-cohomology class. Recall that $H^1(G_{l_{k+1}}, \text{Ad}^0\bar{\rho})$ was two-dimensional over $\mathbb{F}_p$ and so this space contains $p+1$ lines. One of these corresponds to the unramified 1-cohomology class and the restriction of our global 1-cohomology class cannot give this line. Of the others, we do not want it to be null. It seems plausible then that there is a $1 - 1/p$ chance that the restriction of the global 1-cohomology class will be independent of the null 1-cohomology class. Given an infinite number of primes in the Chebotarev class one might perhaps expect one such $l_{k+1}$ to exist.

CORNELL UNIVERSITY, ITHACA, NY
*E-mail address*: ravi@math.cornell.edu## REFERENCES

[B1]   N. BOSTON, Deformation theory of Galois representations, Harvard Ph.D Thesis, 1987.
[B2]   ———, Explicit deformation of Galois Representations, *Invent. Math.* **103** (1991), 181–196.
[BM]   N. BOSTON and B. MAZUR, Explicit universal deformations of Galois representations, in *Algebraic Number Theory*, 1–21, *Adv. Stud. Pure Math.*, 1989.
[Da]   H. DARMON, *The Shimura-Taniyama Conjecture (d'apres Wiles)*, Monographes CICMA Report, 1994-02.
[Fl]   M. FLACH, A finiteness theorem for the symmetric square of an elliptic curve, *Invent. Math.* **109** (1992), 307–327.
[FM]   J.-M. FONTAINE and B. MAZUR, Geometric Galois representations, in *Elliptic Curves, Modular Forms, and Fermat's Last Theorem* (J. Coates and S.-T. Yau, eds.), Ser. Number Theory, I, Internat. Press, Cambridge, MA, 1995.
[H]    K. HABERLAND, *Galois Cohomology of Algebraic Number Fields*, VEB Deutscher Verlag der Wissenschaften, Berlin, 1978.
[LO]   J. C. LAGARIAS and A. M. ODLYZKO, Effective versions of the Chebotarev density theorem, in *Algebraic Number Fields: L-functions and Galois Properties*, 409–464 (*Proc. Sympos.*, Univ. Durham, Durham, 1975), Academic Press, London, 1977.
[M1]   B. MAZUR, Deforming Galois representations, in *Galois Groups over* $\mathbb{Q}$, 385–437 (Berkeley, CA, 1987), *Math. Sci. Res. Inst. Publ.* Springer-Verlag, NY, 1989.
[M2]   ———, An "infinite fern" in the universal deformation space of Galois representations, *Collect. Math.* **48** (1997), 155–193.
[M3]   ———, An introduction to the deformation theory of Galois representations, in *Modular Forms and Fermat's Last Theorem*, 253–311 (Boston, MA, 1995), Springer-Verlag, NY, 1997.